\documentclass [reqno]{amsart}
\usepackage {times,mathptm}
\usepackage {amsmath,amssymb}
\usepackage {color,graphicx}
\usepackage [latin1]{inputenc}

\newcommand {\QQ}{{\mathbb Q}}
\newcommand {\RR}{{\mathbb R}}
\newcommand {\ZZ}{{\mathbb Z}}

\DeclareMathOperator {\ord}{ord}
\DeclareMathOperator {\val}{val}
\DeclareMathOperator {\dist}{dist}
\DeclareMathOperator {\supp}{supp}
\DeclareMathOperator {\Div}{Div}

\DeclareSymbolFont {mysymbols}{OMS}{cmsy}{m}{n}
\DeclareMathSymbol {\calI}{\mathalpha}{mysymbols}{`I}

\DeclareSymbolFont {mylargesymbols}{OMX}{cmex}{m}{n}
\DeclareMathSymbol {\dunion}{\mathop}{mylargesymbols}{"60}

\newcommand {\df}[1]{\textsl {#1}}

\newcommand {\preprint}[2]{preprint \discretionary {#1/}{#2}{#1/#2}}

\renewenvironment {enumerate}%
  {\rule{1mm}{0mm}\begin {oldenumerate}%
    \parskip1ex plus0.5ex \itemsep 0mm \parindent 0mm}%
  {\end {oldenumerate}}

\renewenvironment {itemize}%
  {\rule{1mm}{0mm}\begin {olditemize}%
    \parskip1ex plus0.5ex \itemsep 0mm \parindent 0mm}%
  {\end {olditemize}}

\theoremstyle {plain}
\newtheorem {theorem}{Theorem}[section]
\newtheorem {proposition}[theorem]{Proposition}
\newtheorem {lemma}[theorem]{Lemma}
\newtheorem {corollary}[theorem]{Corollary}
\theoremstyle {definition}
\newtheorem {definition}[theorem]{Definition}
\theoremstyle {remark}
\newtheorem {remark}[theorem]{Remark}
\newtheorem {example}[theorem]{Example}

\begin{document}

  \title [A Riemann-Roch theorem in tropical geometry]{A Riemann-Roch theorem in
  tropical geometry}
\author {Andreas Gathmann and Michael Kerber}
\address {Andreas Gathmann, Fachbereich Mathematik, TU Kaiserslautern, Postfach
  3049, 67653 Kaiserslautern, Germany}
\email {andreas@mathematik.uni-kl.de}
\address {Michael Kerber, Fachbereich Mathematik, TU Kaiserslautern, Postfach
  3049, 67653 Kaiserslautern, Germany}
\email {mkerber@mathematik.uni-kl.de}

\begin {abstract}
  Recently, Baker and Norine have proven a Riemann-Roch theorem for finite
  graphs. We extend their results to metric graphs and thus establish a
  Riemann-Roch theorem for divisors on (abstract) tropical curves.
\end {abstract}

\maketitle

  Tropical algebraic geometry is a recent branch of mathematics that establishes
deep relations between algebro-geometric and purely combinatorial objects.
Ideally, every construction and theorem of algebraic geometry should have a
tropical (i.e.\ combinatorial) counterpart that is then hopefully easier to
understand --- e.g.\ the tropical counterpart of $n$-dimensional varieties are
certain $n$-dimensional polyhedral complexes. In this paper we will establish a
tropical counterpart of the well-known Riemann-Roch theorem for divisors on
curves.

Let us brief\/ly describe the idea of our result. Following Mikhalkin, an
(abstract) tropical curve is simply a connected metric graph $ \Gamma $. A
rational function on $ \Gamma $ is a continuous, piecewise linear real-valued
function $f$ with integer slopes. For such a function and any point $ P \in
\Gamma $ the order $ \ord_P f $ of $f$ in $P$ is the sum of the slopes of $f$
for all edges emanating from $P$. For example, the following picture shows a
rational function $f$ on a tropical curve $ \Gamma $ with simple zeroes at $
P_2 $ and $ P_5 $ (i.e.\ $ \ord_{P_2} f = \ord_{P_5} f = 1 $), and simple poles
at $ P_3 $ and $ P_4 $ (i.e.\ $ \ord_{P_3} f = \ord_{P_4} f = -1 $).

\begin {center} \input {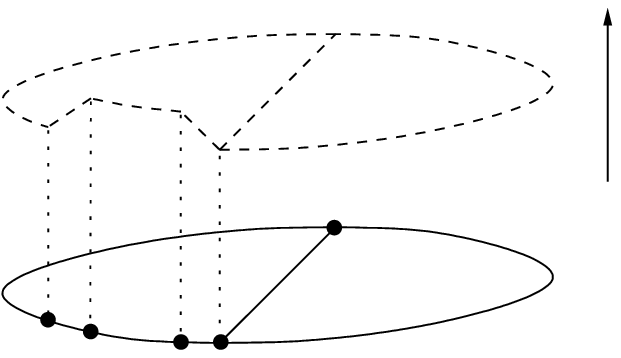} \end {center}

As expected from classical geometry, a divisor on $ \Gamma $ will simply be a
formal $ \ZZ $-linear combination of points of $ \Gamma $. Any rational
function $f$ on $ \Gamma $ gives rise to a divisor $ (f) := \sum_{P \in \Gamma}
\ord_P f \cdot P $ (so that $ (f) = P_2-P_3-P_4+P_5 $ in the above example).

For a given divisor $D$ we denote by $ R(D) $ the space of all rational
functions $f$ on $ \Gamma $ such that $ (f)+D $ is effective, i.e.\ contains
only non-negative coefficients (e.g.\ $ f \in R(P_3+P_4) $ in the example
above). A Riemann-Roch theorem should make a statement about the dimension of
these spaces. However, we will see that in general $ R(D) $ is a polyhedral
complex which is not of pure dimension. As a replacement for the dimension of $
R(D) $ we define $ r(D) $ to be the biggest integer $n$ such that $
R(D-P_1-\cdots-P_n) $ is non-empty for all choices of $ P_1,\dots, P_n \in
\Gamma $ (a number that is closely related to the dimension of the cells of $
R(D) $ as we will see).

With these notations our Riemann-Roch theorem now simply and expectedly states
that
  \[ r(D) - r(K-D) = \deg D + 1 - g, \]
where $g$ is the first Betti number of $ \Gamma $, $ \deg D $ is the degree of
$D$, and $K$ is the canonical divisor of $ \Gamma $ following Zhang \cite {Z},
i.e.\ the sum of all vertices of $ \Gamma $ counted with multiplicity equal to
their respective valence minus 2 (so that $ K=P_1+P_2 $ in our example above).

Our proof relies heavily on a recent result of Baker and Norine that
establishes an analogous result for integer-valued functions on the vertices of
a (non-metric) graph \cite {BN}. Basically, we will interpret this result as a
statement about tropical curves whose edge lengths are integers (so-called
$ \ZZ $-graphs) and rational functions on them whose divisors consist of points
with integer coordinates (so-called $ \ZZ $-divisors). We then pass from
integer to rational and finally real coordinates, as well as to possibly
infinite edge lengths, to establish our Riemann-Roch theorem for tropical
curves.

More precisely, we will first introduce our basic objects of study, namely
divisors and rational functions (and their moduli spaces) on tropical curves in
section \ref {sec-divisors}. We then use the result of Baker and Norine to
prove a Riemann-Roch theorem for $ \ZZ $- and $ \QQ $-divisors in section \ref
{sec-riro-zq} and extend this result in section \ref {sec-riro} to arbitrary
divisors and graphs (with possibly unbounded edges), with the main result being
corollary \ref {riro-final}.

Shortly after this manuscript had appeared on the e-print archive, Mikhalkin
and Zharkov published a preprint that also includes a proof of the Riemann-Roch
theorem for tropical curves \cite {MZ}. Their results have been obtained
independently and without our knowledge, and in fact their method of proof is
entirely different from ours, using Jacobians of tropical curves.

  \section {Tropical rational functions and divisors} \label {sec-divisors}

We start by introducing the basic notations used in this paper, in particular
the notions of (abstract) tropical curves as well as rational functions and
divisors on them.

\begin {definition}[Graphs] \label {def-graph}
  A \df {graph} $\Gamma$ will always mean a finite and connected multigraph,
  not necessarily loop-free (i.e.\ there may be edges that connect a vertex to
  itself). The sets of vertices and edges of $\Gamma$ are denoted $ V(\Gamma) $
  and $ E(\Gamma) $, respectively. The valence of a vertex $ P \in V(\Gamma) $
  will be denoted $ \val (P) $.
  \begin {enumerate}
  \item \label {def-graph-a}
    A \df {metric graph} is a pair $ (\Gamma,l) $ consisting of a graph $
    \Gamma $ together with a length function $ l:E(\Gamma)\rightarrow \RR_{>0}
    $. We identify an edge $e$ with the real interval $ [0,l(e)] $, leading to
    a ``geometric representation'' of the graph by gluing these intervals
    together at their boundary points according to the combinatorics of $
    \Gamma $. By abuse of notation we will usually denote this geometric
    representation also by $ \Gamma $. In this metric space the distance
    between points as well as the distance from a point to a subset will be
    written as $ \dist ( \;\cdot\; , \;\cdot\; ) $. The first Betti number of $
    \Gamma $ will be called the \df {genus} of $ \Gamma $.
  \item \label {def-graph-b}
    If all edge lengths of a metric graph $ \Gamma $ are integers (resp.\
    rational numbers) we call $ \Gamma $ a $ \ZZ $\df {-graph} (resp.\ $ \QQ
    $\df {-graph}). In this case the points of (the geometric representation
    of) $ \Gamma $ with integer (resp.\ rational) distance to the vertices are
    called $ \ZZ $\df {-points} (resp.\ $ \QQ $\df {-points}) of $ \Gamma $.
    We denote the set of these points by $ \Gamma_\ZZ $ and $ \Gamma_\QQ $,
    respectively.
  \item \label {def-graph-c}
    A \df {tropical curve} is a ``metric graph with possibly unbounded ends'',
    i.e.\ a pair $ (\Gamma,l) $ as in \ref {def-graph-a} where the length
    function takes values in $ \RR_{>0} \cup \{\infty\} $, and where each edge
    of length $ \infty $ is identified with the real interval $ [0,\infty] =
    \RR_{\ge 0 } \cup \{ \infty \} $ in such a way that the $ \infty $ end of
    the edge has valence 1. These infinity points of will be called the \df
    {(unbounded) ends} of $ \Gamma $.
  \end {enumerate}
\end {definition}

\begin {remark}
  Note that (in contrast to some other conventions on abstract tropical curves
  found in the literature) our definition allows vertices of valence 1 and 2,
  and adds ``points at infinity'' at each unbounded edge. Note also that every
  metric graph is a tropical curve.
\end {remark}

\begin {definition}[Divisors] \label {def-divisor}
  A \df {divisor} on a tropical curve $\Gamma$ is an element of the free
  abelian group generated by the points of (the geometric representation of) 
  $\Gamma$. The group of all divisors on $\Gamma$ is denoted $ \Div(\Gamma) $.
  The \df {degree} $ \deg D $ of a divisor $ D = \sum_i a_i P_i $ (with $ a_i
  \in \ZZ $ and $ P_i \in \Gamma $) is defined to be the integer $ \sum_i a_i
  $ and obviously gives rise to a morphism $ \deg: \Div (\Gamma) \to \ZZ $. The
  \df {support} $ \supp D $ of $D$ is defined to be the set of all points of $
  \Gamma $ occurring in $D$ with a non-zero coefficient. A divisor is called
  \df {effective} if all its coefficients $ a_i $ are non-negative. On a $ \ZZ
  $-graph (resp.\ $ \QQ $-graph) a divisor $D$ will be called a $ \ZZ $\df
  {-divisor} (resp.\ $ \QQ $\df {-divisor}) if $ \supp D \subset \Gamma_\ZZ $
  (resp.\ $ \supp D \subset \Gamma_\QQ $). Following Zhang \cite {Z} we define
  the canonical divisor of $ \Gamma $ to be
    \[ K_\Gamma := \sum_{P \in V(\Gamma)} (\val (P) - 2) \cdot P; \]
  on a $ \ZZ $-graph (resp.\ $ \QQ $-graph) it is obviously a $ \ZZ $-divisor
  (resp.\ $ \QQ $-divisor).
\end {definition}

\begin {definition}[Rational functions] \label {def-ratfunc}
  A \df {rational function} on a tropical curve $ \Gamma $ is a continuous
  function $ f:\Gamma\rightarrow \RR \cup \{ \pm \infty \} $ such that the
  restriction of $f$ to any edge of $\Gamma$ is a piecewise linear integral
  affine function with a finite number of pieces. In particular, $f$ can take
  on the values $ \pm \infty $ only at the unbounded ends of $ \Gamma $.

  For a rational function $f$ as above and a point $ P \in \Gamma $ the \df
  {order} $ \ord_P f \in \ZZ $ of $f$ at $P$ will be the sum of the outgoing
  slopes of all segments of $ \Gamma $ emanating from $P$ (of which there are $
  \val (P) $ if $ P \in V(\Gamma) $ and 2 otherwise). In particular, if $P$ is
  an unbounded end of $\Gamma$ lying on an unbounded edge $e$ then the order of
  $f$ at $P$ equals the negative of the slope of $f$ at a point on $e$
  sufficiently close to $P$.

  Note that $ \ord_P f = 0 $ for all points $ P \in \Gamma \backslash V(\Gamma)
  $ at which $f$ is locally linear and thus for all but finitely many points.
  We can therefore define the \df {divisor associated to} $f$
    \[ (f) := \sum_{P \in \Gamma} \ord_P f \cdot P \in \Div(\Gamma) \]
  as in classical geometry.
\end {definition}

\begin {remark} \label {deg-zero}
  If $f$ is a rational function on a tropical curve $ \Gamma $ then the degree
  of its associated divisor $ (f) $ is $ \deg (f) = \sum_{P \in \Gamma} \ord_P
  f $. By definition of the order this expression can be written as a sum over
  all segments of $ \Gamma $ on which $f$ is linear, where each such segment
  counts with the sum of the outgoing slopes of $f$ on it at the two end points
  of the segment. But as these two slopes are obviously just opposite numbers
  on each such edge we can conclude that $ \deg (f) = 0 $ --- again analogous
  to the case of compact curves in classical geometry.
\end {remark}

\begin {definition}[Spaces of functions associated to a divisor]
    \label {def-spaces}
  Let $D$ be a divisor of degree $n$ on a tropical curve $ \Gamma $.
  \begin {enumerate}
  \item \label {def-spaces-a}
    We denote by $ R(D) $ the set of all rational functions $f$ on $ \Gamma $
    such that the divisor $ (f)+D $ is effective. Note that for any such $ f
    \in R(D) $ the divisor $ (f)+D $ is a sum of exactly $ \deg ((f)+D) = \deg
    D = n $ points by remark \ref {deg-zero}. So if we define
    \begin {align*}
      S(D) := & \{ (f,P_1,\dots,P_n) ;\;
          \mbox {$f$ a rational function on $ \Gamma $,} \\
        & \qquad \mbox {$ P_1,\dots,P_n \in \Gamma $ such that $ (f)+D =
          P_1+\cdots+P_n $} \}
    \end {align*}
    then we obviously have $ R(D) = S(D) / S_n $, where the symmetric group $
    S_n $ acts on $ S(D) $ by permutation of the points $ P_i $.
  \item \label {def-spaces-b}
    If $ \Gamma $ is a $ \ZZ $-graph and $D$ a $ \ZZ $-divisor we define
    a ``discrete version'' of \ref {def-spaces-a} as follows: let $ \tilde R(D)
    $ be the set of all rational functions $f$ on $ \Gamma $ such that $ (f)+D
    $ is an effective $ \ZZ $-divisor, and set
    \begin {align*}
      \tilde S(D) := & \{ (f,P_1,\dots,P_n) ;\;
          \mbox {$f$ a rational function on $ \Gamma $,} \\
        & \qquad \mbox {$ P_1,\dots,P_n \in \Gamma_\ZZ $ such that $ (f)+D =
          P_1+\cdots+P_n $} \},
    \end {align*}
    so that again $ \tilde R(D) = \tilde S(D)/S_n $.
  \end {enumerate}
  If we want to specify the curve $ \Gamma $ in the notation of these spaces we
  will also write them as $ R_\Gamma(D) $, $ S_\Gamma(D) $, $ \tilde
  R_\Gamma(D) $, and $ \tilde S_\Gamma(D) $, respectively.
\end {definition}

\begin {remark} \label {rem-spaces}
  The spaces $ R(D) $, $ S(D) $, $ \tilde R(D) $, $ \tilde S(D) $ of definition
  \ref {def-spaces} have the following obvious properties:
  \begin {enumerate}
  \item \label {rem-spaces-a}
    all of them are empty if $ \deg D < 0 $;
  \item \label {rem-spaces-b}
    $ R(D-P) \subset R(D) $ and $ \tilde R(D-P) \subset \tilde R(D) $ for all
    $ P \in \Gamma $;
  \item \label {rem-spaces-c}
    $ \tilde R(D) \subset R(D) $ and $ \tilde S(D) \subset S(D) $ if $D$ is a $
    \ZZ $-divisor on a $ \ZZ $-graph $ \Gamma $.
  \end {enumerate}
\end {remark}

We want to see now that $ R(D) $ and $ S(D) $ are \textsl {polyhedral
complexes} in the sense of \cite {GM}, i.e.\ spaces that can be obtained by
gluing finitely many polyhedra along their boundaries, where a polyhedron is
defined to be a subset of a real vector space given by finitely many linear
equalities and strict inequalities. To do this we first need a lemma that
limits the combinatorial possibilities for the elements of $ R(D) $ and $ S(D)
$. For simplicity we will only consider the case of metric graphs here (but it
is in fact easy to see with the same arguments that lemmas \ref {bounded-slope}
and \ref {poly} hold as well for tropical curves, i.e.\ in the presence of
unbounded ends).

\begin {lemma} \label {bounded-slope}
  Let $ p>0 $ be an integer, and let $f$ be a rational function on a metric
  graph $ \Gamma $ that has at most $p$ poles (counted with multiplicities).
  Then the absolute value of the slope of $f$ at any point of $ \Gamma $ (which
  is not a vertex and where $f$ is differentiable) is bounded by a number that
  depends only on $p$ and the non-metric graph $ \Gamma $ (i.e.\ the
  combinatorics of $ \Gamma $).
\end {lemma}

\begin {proof}
  To simplify the notation of this proof we will consider all zeroes and poles
  of $f$ to be vertices of $ \Gamma $ (by making them into 2-valent vertices in
  case they happen to lie in the interior of an edge).

  Let $e$ be any edge of $ \Gamma $ on which $f$ is not constant. Construct a
  path $ \gamma $ along $ \Gamma $ starting with $e$ in the direction in which
  $f$ is increasing, and then successively following the edges of $ \Gamma $,
  at each vertex continuing along an edge on which the outgoing slope of $f$ is
  maximal.

  By our convention on 2-valent vertices above the function $f$ is affine
  linear on each edge of $ \Gamma $. Let us now study how the slope of $f$
  changes along $ \gamma $ when we pass a vertex $ P \in \Gamma $. By
  definition we have $ \lambda_1 + \cdots + \lambda_n = \ord_P f $, where $
  \lambda_1,\dots,\lambda_n $ are the outgoing slopes of $f$ on the edges $
  e_1,\dots,e_n $ adjacent to $P$. Now let $N$ be the maximal valence of a
  vertex occurring in $ \Gamma $, and assume that our path $ \gamma $
  approaches $P$ along the edge $ e_1 $ on which $f$ has incoming slope $
  -\lambda_1 $ greater or equal to $ (N+p)^\alpha $ for some $ \alpha \ge 1 $.
  It then follows that
  \begin {align*}
    \lambda_2 + \cdots + \lambda_n
      &= -\lambda_1 + \ord_P f \\
      &\ge (N+p)^\alpha - p \\
      &= N \, (N+p)^{\alpha-1} + p \, ((N+p)^{\alpha-1} - 1) \\
      &\ge N \, (N+p)^{\alpha-1},
  \end {align*}
  which means that the biggest of the numbers $ \lambda_2,\dots,\lambda_n $,
  i.e.\ the outgoing slope of $f$ along $ \gamma $ when leaving $P$, is at
  least $ (N+p)^{\alpha-1} $ (recall that $ n \le N $ and that $ \lambda_1 $
  can never be the biggest of the $ \lambda_1,\dots,\lambda_n $ since it is
  negative by assumption whereas at least one of the $
  \lambda_2,\dots,\lambda_n $ is positive).

  So if we assume that the slope of $f$ is at least $ (N+p)^\alpha $ on the
  edge $e$ this means by induction that the slope of $f$ on $ \gamma $ is at
  least $ (N+p)^{\alpha-i} $ after crossing $i$ vertices, i.e.\ in particular
  that $f$ is strictly increasing on the first $ \alpha+1 $ edges of $ \gamma
  $. But this is only possible if $ \alpha $ is less than the number of edges
  of $ \Gamma $: otherwise at least one edge must occur twice among the
  first $ \alpha+1 $ edges of $ \gamma $, in contradiction to $f$ being
  strictly increasing on $ \gamma $ in this range. As the initial edge $e$
  was arbitrary this means that the slope of $f$ on any edge is bounded by $
  (N+p)^\alpha $, with $ \alpha $ being the number of edges of $ \Gamma $.
\end {proof}

\begin {lemma} \label {poly}
  For any divisor $D$ on a metric graph $ \Gamma $ the spaces $ R(D) $ and $
  S(D) $ are polyhedral complexes.
\end {lemma}

\begin {proof}
  We will start with $ S(D) $. For each edge $e$ of $ \Gamma $ we choose an
  adjacent vertex that we will call the starting point of $e$. To each element
  $ (f,P_1,\dots,P_n) $ of $ S(D) $ we associate the following discrete data:
  \begin {enumerate}
  \item \label {poly-a}
    the information on which edge or vertex $ P_i $ lies for all $ i = 1,\dots,
    n $;
  \item \label {poly-b}
    the (integer) slope of $f$ on each edge at its starting point;
  \end {enumerate}
  and the following continuous data:
  \begin {enumerate} \setcounter {enumi}{2}
  \item \label {poly-c}
    the distance of each $ P_i $ that lies on an edge from the starting point
    of this edge;
  \item \label {poly-d}
    the value of $f$ at a chosen vertex.
  \end {enumerate}
  These data obviously determine $f$ uniquely: on each edge we know the
  starting slope of $f$ as well as the position and orders of all zeroes and
  poles, so $f$ can be reconstructed on each edge if its starting value on the
  edge is given. As $ \Gamma $ is connected by assumption we can thus
  reconstruct the whole function from the starting value \ref {poly-d}.

  Since there are only finitely many choices for \ref {poly-a} and \ref
  {poly-b} (use lemma \ref {bounded-slope} for \ref {poly-b}), we get a
  stratification of $ S(D) $ with finitely many strata. The data \ref {poly-c}
  and \ref {poly-d} are given by finitely many real variables in each stratum,
  so each stratum is a subset of a real vector space. Finally, the condition on
  the given data to be compatible is given by several linear equalities and
  inequalities (the distances \ref {poly-c} must be positive and less than the
  length of the corresponding edges, and the values of $f$ at the boundary
  points of the edges must be so that we get a well-defined continuous function
  on $ \Gamma $), so that $ S(D) $ is indeed a polyhedral complex.

  The space $ R(D) $ is then simply the quotient of $ S(D) $ by the affine
  linear action of the permutation group of the $ P_1,\dots,P_n $, and hence is
  a polyhedral complex as well.
\end {proof}

\begin {remark}
  For the spaces $ \tilde R(D) $ and $ \tilde S(D) $ the same argument as in
  the proof of lemma \ref {poly} holds, with the only exception that the data
  \ref {poly-c} becomes discrete since the points in $ (f) $ are required to be
  $ \ZZ $-points. Hence the only continuous parameter left is the additive
  constant \ref {poly-d}, i.e.\ both $ \tilde R(D) $ and $ \tilde S(D) $ are
  finite unions of real lines. We can thus regard $ \tilde R(D) $ and $ \tilde
  S(D) $ as ``discrete versions'' of the spaces $ R(D) $ and $ S(D) $.
\end {remark}

The following example shows that the polyhedral complexes $ R(D) $ and $ S(D) $
are in general not pure-dimensional, i.e.\ there may exist inclusion-maximal
cells of different dimensions:

\begin {example} \label {ex-rd-1}
  Consider the canonical divisor $ K_\Gamma=P+Q $ of the metric graph $\Gamma$
  obtained by connecting two cycles $C_1$ and $C_2$ of length 1 by an edge $e$
  of length $ l(e) \in \ZZ_{>0} $ (see the picture below). Furthermore, let $f$
  be a rational function on $\Gamma$ such that $(f)+K_\Gamma = P_1+P_2$.

  Assume first that both $P_1$ and $P_2$ lie in the interior of the edge $e$.
  Note that for all such choices of the points $P_i$ there exists (up to an
  additive constant) exactly one rational function with zeros at $P_1$ and
  $P_2$ and poles at the prescribed points $P$ and $Q$. It follows that the
  corresponding cell in $S(K_\Gamma)$ can be identified with $[0,l(e)] \times
  [0,l(e)] \times \RR $, where the first two factors represent the position of
  the points $P_1$ and $ P_2 $, and the last factor parametrizes the additive
  constant. Hence the dimension of this cell in $ S(K_\Gamma) $ is 3.

  \begin {center} \input {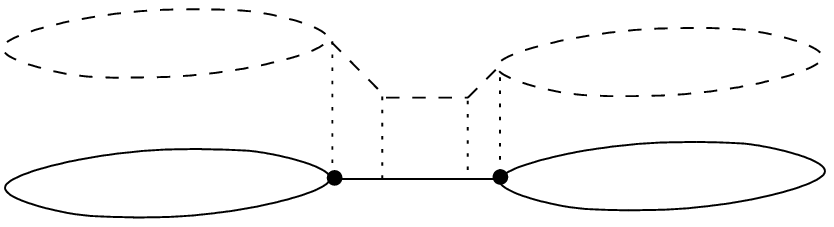} \end {center}

  Next, assume that $ (f,P_1',P_2')\in S(K_\Gamma) $ such that $ P_1' $ is not
  on the closure of $e$ but rather in the interior of a cycle $ C_i $. We will
  see in lemma \ref {not-one-nz} that $ P_2' $ must then lie on the same cycle.
  Moreover, it is easy to check that this requires $ P_2' $ to be the point
  on $ C_i $ ``opposite'' to $ P_1' $ as in the following picture:

  \begin {center} \input {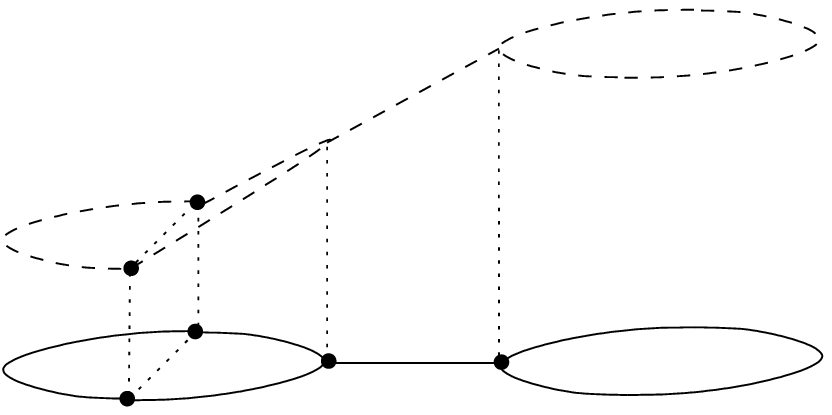} \end {center}

  Hence for each choice of $ P_1' $ on one of the cycles there exists exactly
  one point $ P_2' $ such that $ (f,P_1',P_2') \in S(K_\Gamma) $. It follows
  that this cell of $ S(K_\Gamma) $ can be identified with $ C_i \times \RR $,
  where the second factor parametrizes the additive constant as above. In
  particular, the dimension of this cell is 2.

  Putting all this we obtain the following schematic picture of the polyhedral
  complex $ S(K_\Gamma) $, where for simplicity we have omitted the factor $
  \RR $ corresponding to the additive constant in all cells:

  \begin {center} \input {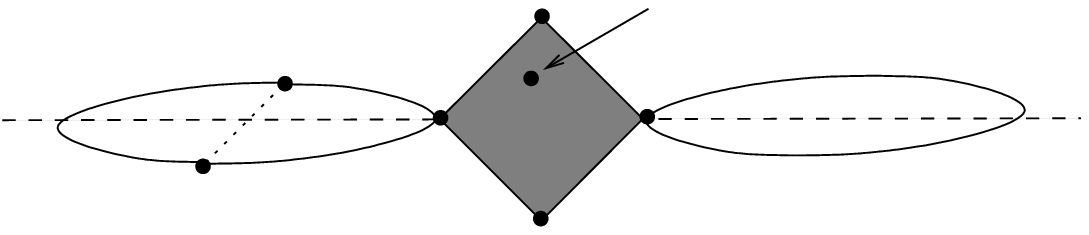} \end {center}

  The space $R(K_\Gamma)$ is then obtained from this by dividing out the action
  of the symmetric group on two elements, which can be realized geometrically
  by ''folding $S(K_\Gamma)$ along the dashed line above''. In particular, both
  $ R(K_\Gamma) $ and $ S(K_\Gamma) $ are not pure-dimensional, but rather have
  components of dimensions 2 and 3.
\end {example}

The above example shows that when formulating a Riemann-Roch type statement
about the dimensions of the spaces $ R(D) $ we have to be careful since these
dimensions are ill-defined in general. The following definition will serve as a
replacement:

\begin {definition} \label {def-rd}
  Let $D$ be a divisor of degree $n$ on a tropical curve $ \Gamma $.
  \begin {enumerate}
  \item \label {def-rd-a}
    We define $ r(D) $ to be the biggest integer $k$ such that for all
    choices of (not necessarily distinct) points $ P_1,\dots,P_k \in \Gamma $
    we have $ R(D-P_1- \cdots -P_k) \neq \emptyset $ (or equivalently $ S
    (D-P_1-\cdots-P_k) \neq \emptyset $), where $ r(D) $ is understood to be
    $ -1 $ if $ R(D) $ (or equivalently $ S(D) $) itself is empty.
  \item \label {def-rd-b}
    If $D$ is a $ \ZZ $-divisor on a $ \ZZ $-graph $ \Gamma $ there is also a
    corresponding ``discrete version'': we let $ \tilde r(D) $ be the biggest
    integer such that $ \tilde R(D-P_1- \cdots -P_k) \neq \emptyset $ for all
    choices of $k$ points $ P_1,\dots,P_k \in \Gamma_\ZZ $.
  \end {enumerate}
  If we want to specify the curve $ \Gamma $ in the notation of these numbers
  we will also write them as $ r_\Gamma(D) $ and $ \tilde r_\Gamma(D) $,
  respectively.
\end {definition}

\begin {example} \label {ex-dim}
  \begin {enumerate}
  \item \label {ex-dim-a}
    By remark \ref {rem-spaces} \ref {rem-spaces-a} it is clear that $ r(D) =
    -1 $ if $ \deg D < 0 $, and $ r(D) \le \deg D $ otherwise. The same
    statement holds for $ \tilde r(D) $ for $ \ZZ $-divisors on $ \ZZ $-graphs.
  \item \label {ex-dim-b}
    For the canonical divisor of the metric graph in example \ref {ex-rd-1} we
    have $ r(K_\Gamma)=1 $ since we have seen that
    \begin {itemize}
    \item
      for all points $ P_1 \in \Gamma $ there is a rational function $f$ with $
      (f)+K_\Gamma = P_1 + P_2 $ (i.e.\ $ f \in R(K_\Gamma-P_1) $);
    \item
      for some choice of $ P_1,P_2 \in \Gamma $ (e.g.\ $ P_1 $ and $ P_2 $ in
      the interior of the circles $ C_1 $ and $ C_2 $, respectively) there is
      no rational function $f$ with $ (f)+K_\Gamma = P_1 + P_2 $.
    \end {itemize}
  \item \label {ex-dim-c}
    Let $ \Gamma $ be a metric graph, and let $ \lambda \in \RR_{>0} $. By a
    \df {rescaling} of $ \Gamma $ by $ \lambda $ we mean the metric graph of
    the same combinatorics as $ \Gamma $ where we replace each edge $e$ of
    length $ l(e) $ by an edge of length $ \lambda \cdot l(e) $. Note that any
    divisor (resp.\ rational function) on $ \Gamma $ gives rise to an induced
    divisor (resp.\ rational function) on the rescaling by also rescaling the
    positions of the points (resp.\ the values of the function). In particular,
    the numbers $ r(D) $ for a divisor $D$ on $ \Gamma $ remain constant under
    rescalings. Note that rescalings by positive integers take $ \ZZ $-graphs
    and $ \ZZ $-divisors again to $ \ZZ $-graphs and $ \ZZ $-divisors, but that
    $ \tilde r(D) $ may change in this case since the rescaling introduces new
    $ \ZZ $-points.
  \end {enumerate}
\end {example}

\begin {remark}
  By the proof of lemma \ref {poly} the continuous parameters for the elements
  $ (f,P_1,\dots,P_n) $ of $ S(D) $ are the positions of the points $ P_i $ and
  the value of $f$ at a chosen vertex. In particular, when passing from $ S(D)
  $ to $ S(D-P) $ for a generic choice of $P$ this fixes one of the $ P_i $ and
  thus makes each cell of $ S(D) $ (disappear or) one dimension smaller. It
  follows that the maximal dimension of the cells of $ S(D) $ (and $ R(D) $) is
  always at least $ r(D)+1 $ (with the $+1$ coming from the additive constant,
  i.e.\ the value of the functions at the chosen vertex).
\end {remark}

\begin {remark} \label {pi-surj}
  There is another interpretation of the numbers $ r(D) $ that we will need
  later: let $D$ be a divisor of degree $n$ on a tropical curve $ \Gamma $,
  let $ i \in \{0,\dots,n\} $, and assume that $ S(D) \neq \emptyset $.
  Consider the forgetful maps
    \[ \pi_i: S(D) \to \Gamma^i ,\quad
         (f,P_1,\dots,P_n) \mapsto (P_1,\dots,P_i). \]
  Note that these maps are morphisms of polyhedral complexes in the sense of
  \cite {GM} (i.e.\ they map each cell of the source to a single cell in the
  target by an affine linear map). It is clear by definition that the number
  $ r(D) $ can be interpreted using these maps as the biggest integer $i$ such
  that $ \pi_i $ is surjective.
\end {remark}

\begin {example}
  Consider again the metric graph $ \Gamma $ of example \ref {ex-rd-1}, but now 
  the spaces $ R(D) $ and $ S(D) $ for the divisor $ D=P'+Q' $, where $ P' $
  are $ Q' $ are interior points of the cycles $ C_1 $ and $ C_2 $,
  respectively. In this case lemma \ref {not-one-nz} will tell us that $
  (f,P_1,P_2) $ can only be in $ S(D) $ if each cycle $ C_i $ contains one of
  the points $ P_1,P_2 $, which is then easily seen to require that in fact $
  \{ P,Q \} = \{P_1,P_2\} $, i.e.\ that $f$ is a constant function. It follows
  that $ R(D) $ is simply the real line, whereas $ S(D) $ is two disjoint
  copies of $ \RR $ (i.e.\ both spaces have pure dimension 1). It also follows
  in the same way that $ r(D)=0 $.

  In particular, when comparing this to the result of examples \ref {ex-rd-1}
  and \ref {ex-dim} \ref {ex-dim-b} (which can be regarded as the limit case
  when $ P' \to P $ and $ Q' \to Q $) we see that $ r(D) $ can jump, and that
  the spaces $ R(D) $ and $ S(D) $ can change quite drastically under
  ``continuous deformations of $D$''. So as in the classical case it is really
  only the number $ r(D) - r(K_\Gamma-D) $, and not $ r(D) $ alone, that will
  turn out to depend on the degree of $D$ and the genus of $ \Gamma $ only.
\end {example}

  \section {Riemann-Roch for $ \QQ $-divisors} \label {sec-riro-zq}

We will now start with the study of Riemann-Roch theorems. Our basic ingredient
is the Riemann-Roch theorem for finite (non-metric) graphs of Baker and Norine
(\cite {BN} theorem 1.11) that is easily translated into our set-up:

\begin {theorem}[Baker and Norine] \label {riro-z}
  Let $ \Gamma $ be a $ \ZZ $-graph of genus $g$ all of whose edge lengths are
  bigger than 1. Then for every $ \ZZ $-divisor $D$ on $ \Gamma $ we have $
  \tilde r(D) - \tilde r(K_\Gamma-D) = \deg D + 1 - g $.
\end {theorem}

\begin {proof}[Sketch of proof]
  We start by replacing each edge $e$ of $ \Gamma $ by a chain of $ l(e) $
  edges of length 1, arriving at a graph whose geometric representation is the
  same as before, and where all $ \ZZ $-points that were in the interior of an
  edge have been turned into 2-valent vertices. Note that by the condition
  that all edge lengths of the original graph are bigger than 1 this implies
  that the new graph has no loops, i.e.\ no edges whose two boundary points
  coincide (an assumption made throughout in \cite {BN}). As it is clear by
  definition that none of the terms in the Riemann-Roch equation changes under
  this transformation it suffices to prove the theorem for the new graph. By
  abuse of notation we will also denote it by $ \Gamma $.

  Note that every rational function $f$ on $ \Gamma $ whose divisor is a $ \ZZ
  $-divisor is uniquely determined by its values on the vertices (since it is
  just given by linear interpolation on the edges). Moreover, up to a possibly
  non-integer global additive constant all these values of $f$ on the vertices
  are integers. Conversely, every integer-valued function on the vertices of $
  \Gamma $ gives rise to a rational function on $ \Gamma $ (by linear
  interpolation) whose divisor is a $ \ZZ $-divisor. As all edge lengths in $
  \Gamma $ are 1 the divisor $ (f) $ can then be rewritten using this
  correspondence as
    \[ (f) = \sum_{\overline {PQ}} \; (f(Q)-f(P)) \cdot (P-Q) \tag {$*$} \]
  where the sum is taken over all edges of $ \Gamma $ (and $P$ and $Q$ denote
  the boundary vertices of these edges in any order). In particular, for a $
  \ZZ $-divisor $D$ the number $ \tilde r(D) $ can also be defined as the
  maximum number $k$ such that for each choice of vertices $ P_1,\dots,P_k $ of
  $ \Gamma $ there is an integer-valued function $f$ on the vertices of $
  \Gamma $ such that $ (f)+D $ is effective, where $(f)$ is defined by $ (*) $.

  This is the approach that Baker and Norine take in \cite {BN}. They establish
  the Riemann-Roch theorem in this set-up, thus proving the theorem as stated
  above. To prove their theorem their first step is to show its equivalence to
  the following two statements:
  \begin {itemize}
  \item $ \tilde{r}(K_\Gamma) \geq g-1 $; and
  \item for any $ \ZZ $-divisor $ D\in\Div(\Gamma) $ there exists a $ \ZZ
    $-divisor $ E\in\Div(\Gamma) $ with $ \deg (E)=g-1 $ and $ \tilde{r}(E)=-1$
    such that exactly one of the sets $\tilde{R}(D)$ and $\tilde{R}(E-D)$ is
    empty.
  \end {itemize}
  The central idea in the proof of these two statements is then to consider
  total orderings on the vertices of $ \Gamma $. For each such ordering there
  is an associated divisor
    \[ E = \sum_{e \in E(\Gamma)} m(e) - \sum_{P \in V(\Gamma)} P \]
  where $ m(e) $ denotes the boundary point of $e$ that is the bigger one in
  the given ordering --- the divisor $E$ in the second statement above can for
  example be taken to be of this form for a suitable ordering (that depends on
  $D$). For details of the proof see \cite {BN}.
\end {proof}

In order to pass from the ``discrete case'' (the spaces $ \tilde R(D) $) to the
``continuous case'' (the spaces $ R(D) $) we need a few lemmas first.

\begin {lemma} \label {not-one-nz}
  Let $D$ be a $\ZZ$-divisor on a $\ZZ$-graph $\Gamma$, and let $
  (f,P_1,\dots,P_n) \in S(D) $. Assume moreover that some $ P_i $ is not a
  $\ZZ$-point. Then on every cycle of $\Gamma$ containing $P_i$ there is
  another point $P_j$ (with $ i \neq j$) that is also not a $\ZZ$-point.
\end {lemma}

\begin {proof}
  Assume that $C$ is a cycle containing exactly one simple zero $ P = P_i \in
  \Gamma \setminus \Gamma_\ZZ $ (note that if $P$ is a multiple zero then we
  are done). Consider the cycle $C$ to be the interval $ [0,l(C)] $ with the
  endpoints identified such that the zero point lies on a vertex, and let $ x
  \in \{1,\dots,l(C)\} $ be the integer such that $ P \in (x-1,x) $ with this
  identification. By adding a suitable constant to $f$ we may assume that $
  f(x-1) \in \ZZ $. Since $ P \in (x-1,x) $ and the slope of $f$ on the
  interval $ [x-1,P] $ differs from that on the interval $ (P,x) $ by $1$ we
  conclude that $ f(x) \notin \ZZ $. As all other points of
  non-differentiability of $f$ on $ [0,l(C)] $ are $ \ZZ $-points by assumption
  it follows that $f(Q) \in \ZZ $ for all $ Q =0,\dots,x-1 $ and $ f(Q) \notin
  \ZZ $ for all $ Q = x,\dots,l(C) $. In particular, we see that $ f(0) \neq
  f(l(C)) $, in contradiction to the continuity of $f$.
\end {proof}

\begin {lemma} \label {lem-move}
  For every $ \ZZ $-divisor $D$ on a $ \ZZ $-graph $ \Gamma $ with $ R(D) \neq
  \emptyset $ we have $ \tilde R(D) \neq \emptyset $.
\end {lemma}

\begin {proof}
  We will prove the statement by induction on $ n := \deg D $.

  Let $ f \in R(D) $, so that $ (f)+D=P_1+\cdots+P_n $ for some (not
  necessarily distinct) points $ P_i \in \Gamma $. In particular, this requires
  of course that $ n \ge 0 $. Moreover, if $ n=0 $ then $ (f)=-D $ is a $ \ZZ
  $-divisor and hence $ f \in \tilde R(D) $. As this finishes the proof in the
  case $ n \le 0 $ we can assume from now on that $ n>0 $, and that the
  statement of the lemma is true for all divisors of degree less than $n$.

  If $ P_i \in \Gamma_\ZZ $ for some $i$ then $ f \in R(D-P_i) $ and hence
  $ \tilde R(D-P_i) \neq \emptyset $ by the induction assumption. As this
  implies $ \tilde R(D) \neq \emptyset $ we have proven the lemma in this case
  and may thus assume from now on that none of the $ P_i $ is a $ \ZZ $-point
  of the curve.

  After possibly relabeling the points $ P_i $ we may assume in addition that $
  0 < \dist (P_n, \Gamma_\ZZ) \le \dist (P_i,\Gamma_\ZZ) $ for all $
  i=1,\dots,n $, i.e.\ that $ P_n $ is a point among the $ P_i $ that minimizes
  the distance to the $ \ZZ $-points of the curve. Let $ P \in \Gamma_\ZZ $ be
  a point with $ \dist (P_n,P) = \dist (P_n,\Gamma_\ZZ) =: d $, and let $
  \Gamma' \subset \Gamma $ be the connected component of $ \Gamma \backslash
  \{P_1,\dots,P_n\} $ that contains $P$. With this notation consider the
  rational function
    \[ h: \Gamma \to \RR, \quad Q \mapsto \begin {cases}
         - \min (d,\dist (Q,\{P_1,\dots,P_n\}))
           & \mbox {if $ Q \in \Gamma' $}, \\
         0 & \mbox {otherwise}.
       \end {cases} \]
  The following picture shows an example of this construction. In this example
  we have assumed for simplicity that all edges of the graph have length 1 so
  that $ \Gamma_\ZZ $ is just the set of vertices. The distance from $ P_5 $ to
  $P$ is smallest among all distances from the $ P_i $ to a vertex, and the
  subset $ \Gamma' \subset \Gamma $ is drawn in bold.
  \begin {center} \input {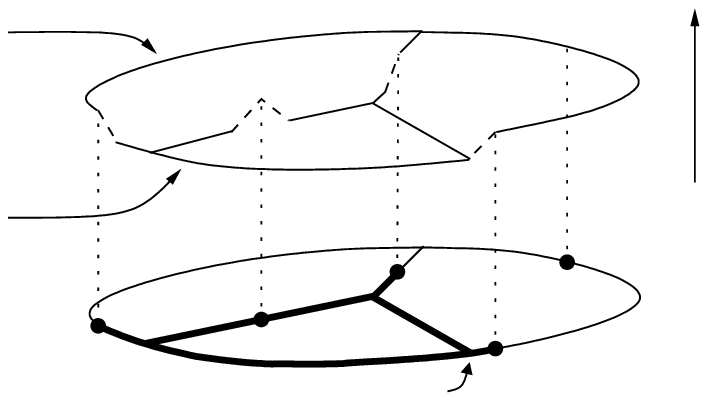} \end {center}
  We claim that $ f+h \in R(D-P) $. In fact, this will prove the lemma since
  $ R(D-P) \neq \emptyset $ implies $ \tilde R(D-P) \neq \emptyset $ and thus
  also $ \tilde R(D) \neq \emptyset $ by the induction assumption.

  To prove that $ f+h \in R(D-P) $ we have to show that $ (f+h)+D-P \ge 0 $,
  or in other words that $ (h)+P_1+\cdots+P_n-P \ge 0 $. Let us assume that
  this statement is false, i.e.\ that there is a point $ Q \in \Gamma $ that is
  contained in the divisor $ (h)+P_1+\cdots+P_n-P $ with a negative
  coefficient. Note that $Q$ cannot be the point $P$ since $ \ord_P h \ge 1 $
  by construction. So $Q$ must be a pole of $h$. But again by construction $h$
  can only have poles at the points $ P_i $, and the order of the poles can be
  at most 2 since the slope of $h$ is $0$ or $ \pm 1 $ everywhere. So the only
  possibility is that $Q$ is a point with $ \ord_Q h = -2 $ that occurs only
  once among the $ P_i $ (as it is the case for $ Q=P_2 $ in the example
  above). But this means that $ \Gamma' $ contains both sides of $Q$, and thus
  (since $ \Gamma' $ is connected) that $ \Gamma' \cup \{Q\} $ contains a cycle
  on which $Q$ is the only point in $ (f) $ that is not a $ \ZZ $-point. But
  this is a contradiction to lemma \ref {not-one-nz} and hence finishes the
  proof of the lemma.
\end {proof}

\begin {proposition} \label {r-tilde-r}
  Let $D$ be a $ \ZZ $-divisor on a $ \ZZ $-graph $ \Gamma $. Then there is an
  integer $ N \ge 1 $ such that $ r(D) = \tilde r(D) $ on every rescaling of $
  \Gamma $ by an integer multiple of $N$ (see example \ref {ex-dim} \ref
  {ex-dim-c}).
\end {proposition}

\begin {proof}
  Let $ n := \deg D $ and $ m := r(D)+1 $, and assume first that $ m \le n $.
  Consider the map $ \pi_m: S(D) \to \Gamma^m $ of remark \ref {pi-surj}.
  As $ \pi_m $ is a morphism of polyhedral complexes its image is closed in $
  \Gamma^m $. Since $ \pi_m $ is not surjective by remark \ref {pi-surj} this
  means that $ \Gamma^m \backslash \pi_m(S(D)) $ is a non-empty open subset of
  $ \Gamma^m $ that consequently must contain an element $ (P_1,\dots,P_m) $
  with rational coordinates. For this element we have $ S(D-P_1-\cdots-P_m) =
  \emptyset $ by construction.
  
  Now let $N$ be the least common multiple of the denominators of these
  coordinates. Then $ P_1,\dots, P_m $ become $ \ZZ $-points on each rescaling
  of $ \Gamma $ by a multiple of $N$, and thus we also have $ \tilde S(D-P_1-
  \cdots-P_m) = \emptyset $ on each such rescaling by remark \ref {rem-spaces}
  \ref {rem-spaces-c}. By definition this then means that $ \tilde r(D) \le m-1
  = r(D) $ on these rescalings. This proves the ``$ \tilde r(D) \le r(D) $''
  part of the proposition in the case $ m \le n $. But note that this part is
  trivial if $ m>n $, since then $ \tilde r(D) \le n \le m-1 = r(D) $ by
  example \ref {ex-dim} \ref {ex-dim-a} (on any rescaling). So we have in fact
  proven the ``$ \tilde r(D) \le r(D) $'' part of the proposition in any case.

  To show the opposite inequality ``$ \tilde r(D) \ge r(D) $'' (which in fact
  holds for any rescaling) we just have to show that $ \tilde R(D-P_1-\cdots
  -P_{r(D)}) \neq \emptyset $ for any choice of $ \ZZ $-points $ P_1,\dots,
  P_{r(D)} $. But this now follows immediately from lemma \ref {lem-move} since
  $ R(D-P_1-\cdots-P_{r(D)}) \neq \emptyset $ by definition.
\end {proof}

We are now ready to prove the Riemann-Roch theorem for $ \QQ $-divisors on $
\QQ $-graphs.

\begin {corollary}[Riemann-Roch for $ \QQ $-graphs] \label {cor-riro-q}
  Let $D$ be a $ \QQ $-divisor on a $ \QQ $-graph $ \Gamma $. Then $ r(D) -
  r(K_\Gamma-D) = \deg D + 1 - g $.
\end {corollary}

\begin {proof}
  Note that it suffices by example \ref {ex-dim} \ref {ex-dim-c} to prove the
  statement after a rescaling of the curve.

  As $ \Gamma $ has only finitely many edges and $D$ contains only finitely
  many points we can assume after such a rescaling that $D$ is in fact a $ \ZZ
  $-divisor on a $ \ZZ $-graph $ \Gamma $, and that all edge lengths of $
  \Gamma $ are bigger than 1. By proposition \ref {r-tilde-r} we can then
  assume after possibly two more rescalings that both $ r(D) = \tilde r(D) $
  and $ r(K_\Gamma-D) = \tilde r(K_\Gamma-D) $. The corollary now follows from
  theorem \ref {riro-z}.
\end {proof}

  \section {Riemann-Roch for tropical curves} \label {sec-riro}

We will now extend our Riemann-Roch theorem for $ \QQ $-graphs (corollary \ref
{cor-riro-q}) in two steps, first to metric graphs (i.e.\ graphs whose edge
lengths need not be rational numbers) and then to tropical curves (i.e.\ graphs
with possibly unbounded edges).

\begin {proposition}[Riemann-Roch for metric graphs] \label {riro-metric}
  For any divisor $D$ on a metric graph $ \Gamma $ of genus $g$ we have
  $ r(D) - r(K_\Gamma - D) = \deg D + 1 - g $.
\end {proposition}

\begin {proof}
  Let $ D = a_1 Q_1 + \cdots + a_m Q_m $, and let $ n = \deg D $. The idea of
  the proof is to find a ``nearby'' $ \QQ $-graph $ \Gamma' $ with a $ \QQ
  $-divisor $D'$ on it such that $ r_{\Gamma'}(D') = r_\Gamma(D) $ and $
  r_{\Gamma'} (K_{\Gamma'}-D') = r (K_\Gamma-D) $, and then to apply the result
  of corollary \ref {cor-riro-q} to this case.

  To do so we will set up a relative version of the spaces $ S(D) $ of
  definition \ref {def-spaces} and the interpretation of $ r(D) $ of remark
  \ref {pi-surj} in terms of these spaces. We fix $ \varepsilon \in \QQ_{>0} $
  smaller than all edge lengths of $ \Gamma $ and denote by $ A(\Gamma) $ the
  set of all metric graphs that are of the same combinatorial type as $ \Gamma
  $ and all of whose edge lengths are greater or equal to $ \varepsilon $. For
  such a metric graph $ \Gamma' \in A(\Gamma) $ we denote by $ B(\Gamma') $ the
  set of all divisors on $ \Gamma' $ that can be written as $ a_1 Q_1' + \cdots
  + a_m Q_m' $ for some $ Q_1',\dots,Q_m' \in \Gamma' $ and the same $
  a_1,\dots,a_m $ as in $D$.
  With these notations we set
  \begin {align*}
    S & := \{ (\Gamma',D',f,P_1,\dots,P_n) ;\;
        \Gamma' \in A(\Gamma), \; D' \in B(\Gamma'), \;
        \mbox {$f$ a rational function on $ \Gamma' $}, \\
      & \qquad \qquad \qquad \qquad
        \mbox {$ P_1,\dots,P_n \in \Gamma' $ such that
        $ (f)+D' = P_1 + \cdots + P_n $} \}, \\
    M_i & := \{ (\Gamma',D',P_1,\dots,P_i) ;\;
        \Gamma' \in A(\Gamma), \; D' \in B(\Gamma'),
        \; P_1,\dots,P_i \in \Gamma' \}
          \qquad \quad \mbox {for $ i=0,\dots,n $}, \\
    M & := \{ (\Gamma',D') ;\;
        \Gamma' \in A(\Gamma), \; D' \in B(\Gamma') \}
  \end {align*}
  In the same way as in lemma \ref {poly} we see that all these spaces are
  polyhedral complexes --- the only difference is that there is some more
  discrete data (corresponding to fixing the edges or vertices on which the
  points in $ D' $ lie) and some more continuous data (corresponding to the
  edge lengths of $ \Gamma' $ and the positions of the points in $ D' $ on
  their respective edges). There are obvious forgetful morphisms of polyhedral
  complexes (i.e.\ continuous maps that send each cell of the source to a
  single cell of the target by an affine linear map)
    \[ \pi_i: S \to M_i, \quad (\Gamma',D',f,P_1,\dots,P_n) \mapsto
                               (\Gamma',D',P_1,\dots,P_i) \]
  and
    \[ p_i: M_i \to M, \quad (\Gamma',D',P_1,\dots,P_i) \mapsto
                             (\Gamma',D'). \]
  As in remark \ref {pi-surj} we have $ r_{\Gamma'}(D') \ge i $ for a divisor $
  D' \in B(\Gamma') $ on a metric graph $ \Gamma' \in A(\Gamma) $ if and only
  if $ \pi_i(S) $ contains $ (\Gamma',D',P_1,\dots,P_i) $ for all $
  P_1,\dots,P_i \in \Gamma' $, or equivalently if and only if $ (\Gamma',D')
  \in M \backslash p_i(M_i \backslash \pi_i(S)) $.

  Since $S$ is a polyhedral complex and $ \pi_i $ a morphism of polyhedral
  complexes it follows that the image $ \pi_i(S) \subset M_i $ is a union
  of closed polyhedra. Consequently, $ M_i \backslash \pi_i(S) $ is a union of
  open polyhedra (i.e.\ an open subset of $ M_i $ whose intersection with each
  polyhedron of $ M_i $ can be written as a union of spaces given by finitely
  many strict linear inequalities).

  Next, note that the map $ p_i $ is open as it is locally just a linear
  projection. It follows that $ p_i (M_i \backslash \pi_i(S)) $, i.e.\ the
  locus in $M$ of all $ (\Gamma',D') $ such that $ r_{\Gamma'}(D')<i $, is
  a union of open polyhedra as well. Consequently, its complement $ M
  \backslash p_i(M_i \backslash \pi_i(S)) $, i.e.\ the locus in $M$ of all $
  (\Gamma',D') $ such that $ r_{\Gamma'}(D') \ge i $, is a union of closed
  polyhedra. Finally, note that all polyhedral complexes and morphisms involved
  in our construction are defined over $ \QQ $, so that the locus of all $
  (\Gamma',D') $ with $ r_{\Gamma'}(D') <i $ (resp.\ $ r_{\Gamma'}(D') \ge i $)
  is in fact a union of \textsl {rational} open (resp.\ closed) polyhedra in
  $M$. Of course, the same arguments hold for $ r_{\Gamma'}(K_{\Gamma'}-D') $
  as well.

  We are now ready to finish the proof of the proposition. By what we have said
  above the locus of all $ (\Gamma',D') $ in $M$ such that $ r_{\Gamma'}(D') <
  r_\Gamma(D)+1 $ and $ r_{\Gamma'}(K_{\Gamma'}-D') < r_\Gamma(K_\Gamma-D) + 1
  $ is an open neighborhood $U$ of $ (\Gamma,D) $. Conversely, the locus of all
  $ (\Gamma',D') $ in $M$ such that $ r_{\Gamma'}(D') \ge r_\Gamma(D) $ and $
  r_{\Gamma'}(K_{\Gamma'}-D') \ge r_\Gamma(K_\Gamma-D) $ is a union $V$ of
  rational closed polyhedra. In particular, this means that the rational points
  of $V$ are dense in $V$. As $ U \cap V $ is non-empty (it contains the point
  $ (\Gamma,D) $) it follows that there is a rational point in $ U \cap V $,
  i.e.\ a $ \QQ $-graph $ \Gamma' $ with a $ \QQ $-divisor $ D' $ on it such
  that $ r_{\Gamma'}(D') = r_\Gamma(D) $ and $ r_{\Gamma'}(K_{\Gamma'}-D') =
  r_\Gamma(K_\Gamma-D) $. As $ \Gamma' $ and $ \Gamma $ have the same genus,
  and $ D' $ and $D$ the same degree, the proposition now follows from
  corollary \ref {cor-riro-q}.
\end {proof}

So far we have only considered metric graphs, i.e.\ tropical curves in which
every edge is of finite length. In our final step of the proof of the
Riemann-Roch theorem we will now extend this result to arbitrary tropical
curves (with possibly infinite edges). In order to do this we will first
introduce the notion of equivalence of divisors.

\begin {definition} \label {def-equiv}
  Two divisors $D$ and $D'$ on a tropical curve $\Gamma$ are called \df
  {equivalent} (written $ D \sim D' $) if there exists a rational function $f$
  on $\Gamma$ such that $ D' = D + (f) $.
\end {definition}

\begin {remark} \label {rem-equiv}
  If $ D \sim D' $, i.e.\ $ D' = D + (f) $ for a rational function $f$, then
  it is obvious that the map $ R(D') \to R(D) $, $ g \mapsto g+f $ is a
  bijection. In particular, this means that $ r(D) = r(D') $, i.e.\ that the
  function $ r: \Div(\Gamma) \to \ZZ $ depends only on the equivalence class of
  $D$.
\end {remark}

\begin {lemma} \label {lem-equiv}
  Let $ \bar \Gamma $ be a tropical curve, and let $ \Gamma $ be the metric
  graph obtained from $ \bar \Gamma $ by removing all unbounded edges. Then
  every divisor $ D \in \Div (\bar \Gamma) $ is equivalent on $ \bar \Gamma $
  to a divisor $ D' $ with $ \supp D' \subset \Gamma $. Moreover, if $D$ is
  effective then $ D' $ can be chosen to be effective as well.
\end {lemma}

\begin {proof}
  To any $ P \in \bar\Gamma $, we associate a rational function $ f_P $ as
  follows. If $ P\in \Gamma $, we define $f_P$ to be the zero function.
  Otherwise, if $P$ lies on some unbounded edge $E$, we define
    \[ f_P: \bar \Gamma \to \RR \cup \{ \infty \}, \quad
         Q \mapsto \begin{cases}
           \min (\dist (P,\Gamma),\dist (Q,\Gamma)) & \text{ if } Q\in E,\\
           0 & \text{ if } Q\not\in E.
                   \end{cases}  \]
  If $ P \notin\Gamma $, then the function $f_P$ has a simple pole at $P$
  and no other zeros or poles away from $ \Gamma $. The following picture shows
  an example of such a function, where the metric graph $ \Gamma $ is drawn in
  bold:

  \begin {center} \input {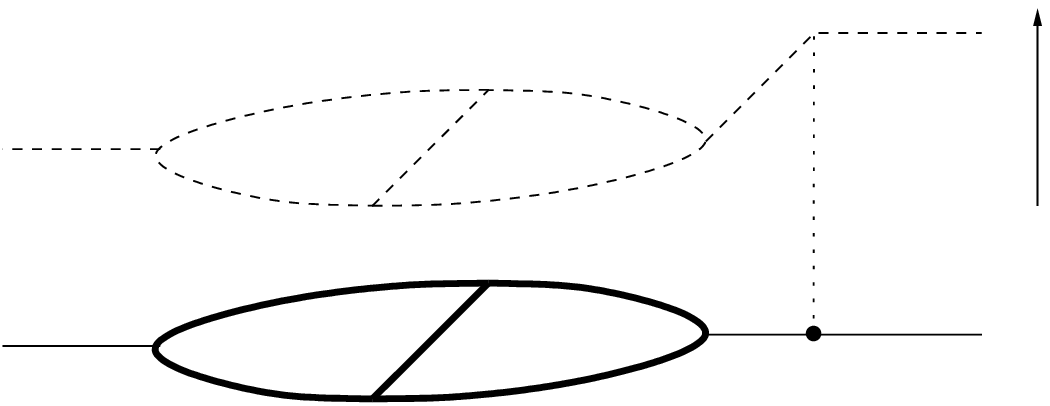} \end {center}

  So if $ D = a_1 P_1 + \cdots + a_n P_n $ and we set $ f = \sum_i a_i
  f_{P_i} $ then $ D+(f) $ is a divisor equivalent to $D$ with no zeros or
  poles away from $ \Gamma $. Moreover, if $D$ is effective then $ D+(f) $ is
  effective as well since all poles of $f$ are cancelled by $D$ by
  construction.
\end{proof}

\begin {remark}\label {KGamma}
  With notations as above, let $ P_1,\cdots,P_n $ denote the end points
  of the unbounded edges $E_i$ of $ \bar\Gamma $, and consider the function
  $f=\sum_i f_{P_i}$. Then $f$ is zero on the graph $\Gamma$ and has slope one
  on each unbounded edge. If we denote for all $ i\in\{1,\dots,n\} $ the point
  $ E_i\cap \Gamma$ by $Q_i$, then $ (f)=\sum Q_i -\sum P_i$. Hence
  $ K_\Gamma+(f)=K_{\bar\Gamma}$, i.e.\ $ K_\Gamma\sim K_{\bar\Gamma}$ on
  $\bar\Gamma$.
\end {remark}

\begin {lemma} \label {lem-remove}
  As in the previous lemma let $ \bar \Gamma $ be a tropical curve, and let $
  \Gamma $ be the metric graph obtained from $ \bar \Gamma $ by removing all
  unbounded edges. Moreover, let $D$ be a divisor on $ \Gamma $ (that can then
  also be thought of as a divisor on $ \bar \Gamma $ with support on $ \Gamma
  $). Then $ R_\Gamma(D) \neq \emptyset $ if and only if $ R_{\Gamma'}(D) \neq
  \emptyset $.
\end {lemma}

\begin {proof}
  ``$ \Rightarrow $'': Let $f$ be a rational function in $R_\Gamma(D)$. Extend
  $f$ to a rational function $ \bar f $ on $ \bar\Gamma $ so that it is
  constant on each unbounded edge. Then $ \bar f \in R_{\bar\Gamma}(D) $.

  ``$ \Leftarrow $'': Let $ \bar f \in R_{\bar\Gamma}(D) $, and set $ f =
  \bar f|_\Gamma $. Let $e$ be an unbounded edge of $ \bar \Gamma $, and let $
  P = \Gamma \cap e $ be the vertex where $e$ is attached to $ \Gamma $. Since
  $ \bar f $ has no poles on $e$ it follows that $ \bar f|_e $ is (not
  necessarily strictly) decreasing if we identify $e$ with the real interval $
  [0,\infty] $. Hence the order of $f$ on $ \Gamma $ at $P$ cannot be less than
  the order of $ \bar f $ on $ \bar \Gamma $ at $P$, and so it follows that $ f
  \in R_\Gamma(D) $.
\end {proof}

\begin {remark} \label {rem-remove}
  Let $ \bar \Gamma $, $ \Gamma $, and $D$ as in lemma \ref {lem-remove}.
  By lemma \ref {lem-equiv} any effective divisor $ P_1 + \cdots + P_k $ on $
  \bar \Gamma $ is equivalent to an effective divisor $ P_1' + \cdots + P_k' $
  with support on $ \Gamma $. So by remark \ref {rem-equiv} the number $
  r_{\bar \Gamma} (D) $ can also be thought of as the biggest integer $k$ such
  that $ R_{\bar \Gamma}(D-P_1-\cdots-P_k) \neq \emptyset $ for all $
  P_1,\dots,P_k \in \Gamma $ (instead of for all $ P_1,\dots,P_k \in \bar
  \Gamma $). By lemma \ref {lem-remove} we can therefore conclude that $
  r_{\bar \Gamma}(D) = r_\Gamma(D) $.
\end {remark}

With these results we are now able to prove our main theorem:

\begin {corollary}[Riemann-Roch for tropical curves]
    \label {riro-final}
  For any divisor $D$ on a tropical curve $ \bar\Gamma $ of genus $g$ we have
  $ r(D) - r(K_{\bar\Gamma} - D) = \deg D + 1 - g $.
\end {corollary}

\begin {proof}
  Let $ \Gamma $ be the metric graph obtained from $ \bar\Gamma $ by removing
  all unbounded edges. By lemma \ref {lem-equiv} and remark \ref {rem-equiv} we
  may assume that $ \supp D \subset \Gamma $. Moreover, by remark
  \ref {KGamma} we can replace $ K_{\bar\Gamma} $ by $ K_\Gamma $ (which
  also has support in $ \Gamma $) in the Riemann-Roch equation. Finally, remark
  \ref {rem-remove} now tells us that we may replace $ r_{\bar \Gamma} (D) $
  and $ r_{\bar \Gamma} (K_\Gamma-D) $ by $ r_\Gamma (D) $ and $ r_\Gamma
  (K_\Gamma-D) $ respectively, so that the statement follows from proposition
  \ref {riro-metric}.
\end {proof}

\end {document}